\pdfoutput=1
\documentclass[12pt]{article}

\usepackage{amsfonts,amsmath,amssymb}
\usepackage{mathrsfs,mathtools,url}
\usepackage{latexsym}
\usepackage{tikz,color}
\usepackage{cite}

\newtheorem{theorem}{Theorem}[section]

\newtheorem{lemma}[theorem]{Lemma}
\newtheorem{proposition}[theorem]{Proposition}

\DeclareMathOperator {\gp} {gp}
\DeclareMathOperator {\diam} {diam}

\let\deg\relax
\DeclareMathOperator {\deg} {deg}

\newcommand{\gpe}{\gp_{\rm e}}
\newcommand{\proof}{\noindent{\bf Proof.\ }}
\newcommand{\qed}{\hfill $\square$ \bigskip}

\vspace{8mm}
\title{Extremal edge general position sets in some graphs}

\author{Jing Tian$\/^{a}$ \and Sandi Klav\v{z}ar$^{b, c, d}$ \and Elif Tan$^{e}$ \and    \\\\
$^{a}$ \small  School of Science, Zhejiang University of Science and Technology, \\
\small Hangzhou, Zhejiang 310023, PR China\\
\small {\tt jingtian526@126.com}\\
$^{b}$ \small Faculty of Mathematics and Physics, University of Ljubljana, Slovenia\\
\small {\tt sandi.klavzar@fmf.uni-lj.si}\\
$^{c}$ \small Faculty of Natural Sciences and Mathematics, University of Maribor, Slovenia\\
$^{d}$ \small Institute of Mathematics, Physics and Mechanics, Ljubljana, Slovenia \\
$^e$ \small Department of Mathematics, Ankara University, Ankara, Turkey  \\
\small {\tt etan@ankara.edu.tr}\\
}
\date{}

\begin{document}

\maketitle

\begin{abstract}
A set of edges $X\subseteq E(G)$ of a graph $G$ is an edge general position set if no three edges from $X$ lie on a common shortest path. The edge general position number ${\rm gp}_{\rm e}(G)$ of $G$ is the cardinality of a largest edge general position set in $G$. Graphs $G$ with ${\rm gp}_{\rm e}(G) = |E(G)| - 1$ and with ${\rm gp}_{\rm e}(G) = 3$ are respectively characterized. Sharp upper and lower bounds on ${\rm gp}_{\rm e}(G)$ are proved for block graphs $G$ and exact values are determined for several specific block graphs.
\end{abstract}

\noindent
\textbf{Keywords:} general position set; edge general position set; cut-vertex; diametral path; block graph

\medskip\noindent
\textbf{AMS Math.\ Subj.\ Class.\ (2020)}: 05C12, 05C35


\section{Introduction}

Edge general position sets in graphs occur as the edge variant of general position sets. As we know, general position sets have already been extensively researched, see the seminal papers~\cite{manuel-2018, ullas-2016}, some of the subsequent ones~\cite{klavzar-2023a, klavzar-2021, korze-2023, patkos-2020, tian-2021a, tian-2023+, yao-2022}, and references therein. The edge version has been introduced in~\cite{manuel-2022} and further studied in~\cite{klavzar-2023b, manuel-2023+}. In this paper we continue this direction of research.

Let $G = (V(G), E(G))$ be a graph. Then we say that $X\subseteq E(G)$ is an {\em edge general position set} if for any shortest path $P$ of $G$ and any three edges $e, e', e''\in X$ we have $|E(P) \cap \{e, e', e''\}| \le 2$.  The {\em edge general position number}  $\gpe(G)$ of $G$ is the cardinality of a largest edge general position set. An edge general position set of cardinality $\gpe(G)$ is called a {\em $\gpe$-set} of $G$.

The main results of the seminal paper~\cite{manuel-2022} are the formulas for the  edge general position number of hypercubes and of Cartesian grids. As a point of interest, we would like to add that the general position number of hypercubes is a notorious difficult problem, cf.~\cite{korner-1995}. The paper is organized as follows: in Section~\ref{sec:extreme}, we characterize the graphs $G$ with $\gpe(G) = |E(G)| - 1$ and the graphs $G$ with $\gpe(G) = 3$. We also discuss the graphs $G$ with $\gpe(G) = 4$. In Section~\ref{sec:block} we focus on block graphs. We prove sharp upper and lower bounds and determine exact values in several specific cases. In the rest of the introduction we list some definitions and state a couple of results needed later on.

Unless stated otherwise, the graphs considered are connected. If $G$ is a graph and $R\subseteq V(G)$, then the subgraph of $G$ induced by $R$ is denoted by $G[R]$. The {\em distance} $d_G(u,v)$ between vertices $u$ and $v$ of $G$ is the number of edges on a shortest $u,v$-path. A subgraph $H$ of $G$ is {\em isometric} if for each pair of vertices $u,v\in V(H)$ we have $d_H(u,v) = d_G(u,v)$. The {\em diameter} of $G$ is  the maximum distance between pairs of vertices of $G$ and is denoted by $\diam(G)$. A shortest $u,v$-path in $G$ is {\em diametral} if $d_G(u,v) = \diam(G)$. The order, the size, and the maximum degree of $G$ are denoted by $n(G)$, $m(G)$, and $\Delta(G)$, respectively. A \textit{block} of $G$ is a maximal connected subgraph of $G$ that has no cut-vertex.

If $G$ is a graph, then the edges incident to a vertex form an edge general position set, hence we have
\begin{equation}
\gpe(G) \ge \Delta(G)\,.
\label{eq:Delta}
\end{equation}
For later use we also recall the following results from~\cite{manuel-2022}. For its statement recall that an edge of a tree is a {\em pendant edge} if it is incident with a leaf and that $C_n$, $n\ge 3$, denotes the cycle of order $n$.  

\begin{proposition}
\label{pro:earlier-results}
(i) $\gpe(C_n)=n$ for $3\leq n\leq 5$ and $\gpe(C_n)=4$ for $n\geq 6$. \\
(ii) If $L$ is the set of pendant edges of a tree $T$, then $\gpe(T)=|L|$.
\end{proposition}

\section{Graphs with extremal edge general position numbers}
\label{sec:extreme}

In the main result of this section we characterize the graphs $G$ with $\gpe(G) = m(G) - 1$. We further characterize the graphs with the edge general position number equal to $3$ and show that the variety of graphs with this number equal to $4$ is large.

In~\cite{manuel-2022} it was observed that if $\diam(G) = 2$, then $\gpe(G) = m(G)$. A slightly more general formulation of this fact reads as follows.

\begin{lemma}
\label{lem:diam2}
If $G$ is a graph, then $\gpe(G) = m(G)$ if and only if $\diam(G) \le 2$.
\end{lemma}

To characterize graphs $G$ with  $\gpe(G) = m(G) - 1$, we introduce the following two families of graphs.

The family ${\cal G}_1$ consists of all the graphs $G$ that can be obtained from an arbitrary graph $H$ with $\diam(H) = 2$ by attaching a pendant edge to a vertex $u$ of $H$ with $\deg_H(u) \le n(H) - 2$. We will say that this pendant edge of $G$ is the {\em special edge} of $G$. Knowing how the graph $G$ from ${\cal G}_1$ was constructed, there is no doubt which of its edges is the special edge. On the other hand, if only $G$ is given, then it might not be clear which of its edges is the special one. For instance, $K_4$ belongs to
${\cal G}_1$ since it can be obtained from $P_3$ by attaching a pendant edge to its vertex of degree $1$. Then each of the two pendant edges of $P_4$ can be declared as the special one. In such a case we select one such edge and fix it, so that the special edge of $G$ is well-defined. In Fig.~\ref{fig:G1-and-G2} the graph $Z_1$ belongs to the family ${\cal G}_1$, where the top edge is the special edge of $Z_1$.

The family ${\cal G}_2$ consists of all the graphs $G$ constructed as follows. Let $G_0$, $G_1$, and $G_2$ be arbitrary graphs not necessarily connected. Then $G$ is obtained from the disjoint union of $G_0$, $G_1$, $G_2$, and $K_2$, where $V(K_2) = \{x_1, x_2\}$, by adding all possible edges between $K_2$ and $G_0$, adding the edges between $x_1$ and all vertices from $G_1$, and adding the edges between $x_2$ and all vertices of $G_2$. We will say that the edge $x_1x_2$ of $G$ is the {\em central edge} of $G$. In Fig.~\ref{fig:G1-and-G2} the graph $Z_2$ belongs to the family ${\cal G}_2$. Note that for $Z_2$ we have $G_0 = P_3\cup P_2$, $G_1 = C_4$, and $G_2 = K_2 \cup 2K_1$. Another graph from the family ${\cal G}_2$ is the graph $G'$ from Fig.~\ref{fig:G'-and-G''}, where $G_0$ is the empty graph, $G_1 = K_3$, and $G_2 = K_2$.

\begin{figure}[ht!]
\begin{center}
\begin{tikzpicture}[scale=1.0,style=thick]
\tikzstyle{every node}=[draw=none,fill=none]
\def\vr{3pt} 

\begin{scope}[yshift = 0cm, xshift = 0cm]
    \node [below=0.5mm] at (0.25,0) {};
    \node [below=0.5mm] at (0.5,0.5) {};
    \node [below=0.5mm] at (1.5,0.5) {};
    \node [below=0.5mm] at (1.75,0) {};
    \node [below=0.5mm] at (2.5,1.5) {};
    \node [below=0.5mm] at (1.75,1.5) {};
    \node [below=0.5mm] at (0.25,1.5) {};
    \node [below=0.5mm] at (-0.5,1.5) {};
    \node [below=0.5mm] at (1,2) {};
    \node [below=0.5mm] at (1,2.75) {};
    \node [below=0.5mm] at (1,3.75) {};

\path (0.25,0) coordinate (x1);
\path (1.75,0) coordinate (x2);
\path (2.5,1.5) coordinate (x3);
\path (1,2.75) coordinate (x4);
\path (-0.5,1.5) coordinate (x5);
\path (0.5,0.5) coordinate (x6);
\path (1.5,0.5) coordinate (x7);
\path (1.75,1.5) coordinate (x8);
\path (1,2) coordinate (x9);
\path (0.25,1.5) coordinate (x10);
\path (1,3.75) coordinate (x11);

\draw (x1)--(x2)--(x3) -- (x4)--(x5)--(x1);
\draw (x6)--(x8)--(x10)--(x7)--(x9)--(x6);
\draw (x1)--(x6);
\draw (x2)--(x7);
\draw (x3)--(x8);
\draw (x4)--(x9);
\draw (x5)--(x10);
\draw (x4)--(x11);

\draw (x1)  [fill=white] circle (\vr);
\draw (x2)  [fill=white] circle (\vr);
\draw (x3)  [fill=white] circle (\vr);
\draw (x4)  [fill=white] circle (\vr);
\draw (x5)  [fill=white] circle (\vr);
\draw (x6)  [fill=white] circle (\vr);
\draw (x7)  [fill=white] circle (\vr);
\draw (x8)  [fill=white] circle (\vr);
\draw (x9)  [fill=white] circle (\vr);
\draw (x10)  [fill=white] circle (\vr);
\draw (x11)  [fill=white] circle (\vr);

\draw (1.0,-1) node {$Z_1$};
\end{scope}
\end{tikzpicture}
\hspace{10mm}
\begin{tikzpicture}[scale=1.2,style=thick]
\tikzstyle{every node}=[draw=none,fill=none]
\def\vr{3pt} 

\begin{scope}[yshift = 0cm, xshift = 0cm]
 \node [below=0.5mm] at (0,0.75) {};
 \node [below=0.5mm] at (0,1.5) {};
 \node [below=0.5mm] at (1,2) {};
 \node [below=0.5mm] at (2,1) {};
 \node [below=0.5mm] at (2,1.5) {};
 \node [below=0.5mm] at (2,2.5) {};
 \node [below=0.5mm] at (2,3) {};
 \node [below=0.5mm] at (0.75,3) {};
 \node [below=0.5mm] at (0,3) {};
 \node [below=0.5mm] at (-0.75,3) {};
 \node [below=0.5mm] at (-1,2) {};
 \node [below=0.5mm] at (-2,3) {};
 \node [below=0.5mm] at (-2,2.3) {};
 \node [below=0.5mm] at (-2,1.7) {};
 \node [below=0.5mm] at (-2,1) {};

\path (0,0.75) coordinate (x1);
\path (1,2) coordinate (x2);
\path (2,1) coordinate (x3);
\path (2,1.5) coordinate (x4);
\path (2,2.5) coordinate (x5);
\path (2,3) coordinate (x6);
\path (0.75,3) coordinate (x7);
\path (0,3) coordinate (x8);
\path (-0.75,3) coordinate (x9);
\path (-1,2) coordinate (x10);
\path (-2,3) coordinate (x11);
\path (-2,2.3) coordinate (x12);
\path (-2,1.7) coordinate (x13);
\path (-2,1) coordinate (x14);
\path (0,1.5) coordinate (x15);

\draw (x1)-- (x2)--(x3);
\draw(x2) -- (x4)--(x5)-- (x2)--(x6);
\draw(x2) -- (x7)--(x8)-- (x9)--(x10)--(x11)--(x12)-- (x13)--(x14)--(x10);
\draw (x10)--(x9)-- (x2);
\draw (x10)--(x8)-- (x2);
\draw (x10)--(x7)-- (x2);
\draw (x10)--(x15)-- (x2);
\draw (x10)--(x1)-- (x2);
\draw (x1)--(x15);
\draw (x10)-- (x2);
\draw (x10)--(x11);
\draw (x10)--(x12);
\draw (x10)--(x13);
\draw (x10)--(x14);
\draw (x11) .. controls (-2.5,2).. (x14);

\draw (x1)  [fill=white] circle (\vr);
\draw (x2)  [fill=white] circle (\vr);
\draw (x3)  [fill=white] circle (\vr);
\draw (x4)  [fill=white] circle (\vr);
\draw (x5)  [fill=white] circle (\vr);
\draw (x6)  [fill=white] circle (\vr);
\draw (x7)  [fill=white] circle (\vr);
\draw (x8)  [fill=white] circle (\vr);
\draw (x9)  [fill=white] circle (\vr);
\draw (x10)  [fill=white] circle (\vr);
\draw (x11)  [fill=white] circle (\vr);
\draw (x12)  [fill=white] circle (\vr);
\draw (x13)  [fill=white] circle (\vr);
\draw (x14)  [fill=white] circle (\vr);
\draw (x15)  [fill=white] circle (\vr);

\draw (0,0) node {$Z_2$};
\draw (-0.95,1.65) node {$x_1$};
\draw (1.05,1.65) node {$x_2$};

\end{scope}
\end{tikzpicture}
\end{center}
\caption{Graphs $Z_1\in {\cal G}_1$ and $Z_2\in {\cal G}_2$.}
\label{fig:G1-and-G2}
\end{figure}

\begin{theorem}
\label{thm:m(G)-1}
Let $G$ be a (connected) graph with $n(G) \ge 4$. Then $\gpe(G) = m(G) - 1$ if and only if $G\in {\cal G}_1 \cup {\cal G}_2$.
\end{theorem}

\proof
Assume first that $G\in {\cal G}_1 \cup {\cal G}_2$. Then $\diam(G) = 3$, and hence $\gpe(G) \le m(G) - 1$. If $G\in {\cal G}_1$, then it is clear that $E(G)\setminus \{e\}$ is an edge general position set of $G$, where $e$ is the special edge of $G$. Let next $G\in {\cal G}_2$. Then we claim that $E(G)\setminus \{x_1x_2\}$ is an edge general position set of $G$, where $x_1x_2$ is the central edge of $G$. Indeed, this follows by the fact that every shortest path of $G$ of length $3$ passes through $x_1x_2$. In both cases we have $\gpe(G) \ge m(G) - 1$ and we may conclude that $\gpe(G) = m(G) - 1$.

Conversely, assume that $\gpe(G) = m(G) - 1$.

\medskip\noindent
{\bf Claim 1}: $\diam(G) = 3$.\\
By Lemma~\ref{lem:diam2} we have $\diam(G)\ge 3$. Moreover, if $\diam(G) \ge 4$, then at least two edges of an arbitrary diametral path of $G$ do not lie in an arbitrary edge general position set, hence $\gpe(G) \le m(G)-2$. We conclude that $\diam(G) = 3$.

\medskip\noindent
{\bf Claim 2}: There exists an edge $e\in E(G)$ which lies on every diametral path of $G$.\\
Let $R_1, \ldots, R_\ell$ be the diametral paths of $G$. If $V(R_i) \cap V(R_j) = \emptyset$, then every edge general position set of $G$ misses at least one edge of $R_i$ and at least one edge of $R_j$. Hence each two diametral paths share at least one edge.

Consider now arbitrary three diametral paths $R_i$, $R_j$ and $R_k$. We claim that they have at least one common edge. Suppose on the contrary that this is not the case and let $X$ be a largest edge general position set of $G$. Assume that $E(R_i) \cap E(R_j) = \{e\}$, $E(R_i) \cap E(R_k) = \{f\}$, and $ E(R_j) \cap E(R_k) = \{g\}$. If  $|\{e,f,g\}| \le 2$, where we without loss of generality have $f=g$, then the edge $f$ lies on all three paths. Hence we must have $|\{e,f,g\}| = 3$. Recall that at least one edge of $R_i$ (as well as of $R_j$ and of $R_k$) is not in $X$. If this edge is $e$, then at least one additional edge of $R_j$ is also not in $X$. Similarly, if $f\notin X$, then at least one additional edge of $R_k$ is also not in $X$. Finally, if the unique edge from $V(R_i)\setminus \{f,g\}$ is not in $X$, then on each of $R_j$ and $R_k$ one additional edge is not in $X$. Hence in all the cases we have $\gpe(G) < m(G) - 1$, a contradiction. We conclude that $E(R_i) \cap E(R_j) \cap E(R_k) \ne \emptyset$.

Let $e\in E(R_i) \cap E(R_j) \cap E(R_k)$ and let $R_{k'}$ be an arbitrary additional diametral path. We are going to show that all these four paths have a common edge. If $e\in E(R_{k'})$, there is nothing to prove. Assume hence that $e\notin E(R_{k'})$. We now consider two cases. If $e\in X$, then on each of the paths $R_i$, $R_j$, $R_k$ and $R_{k'}$ there lies at least one edge which is not in $X$. If this is the same edge, we are done. Otherwise there exist two different edges not in $X$, which is not possible as we have assumed that $\gpe(G) = m(G) - 1$. In the second case $e\notin X$. Since $e\notin E(R_{k'})$, the path $R_{k'}$ contains an edge different from $e$ not in $X$, which is again not possible. By induction we now conclude the truth of Claim 2.

\medskip
By Claim~2, $G$ has an edge $e = uv$ which lies on all diametral paths of $G$.

Assume first that $e$ is the first edge of some diametral path with the first vertex of it being $u$.  Then $e$ must be the first edge of every diametral path, for otherwise we would have $\diam(G) \ge 4$ or there would exist another diametral path which would not contain $e$. Moreover, by the same reason we see that $\deg_G(u) = 1$. Let $H = G - u$. Then $\diam(H) = 2$ and hence, since $\diam(G) = 3$, it follows that $\deg_H(v) \le n(H) - 2$. We conclude that $G\in {\cal G}_1$.

Assume second that $e=uv$ is the middle edge of some diametral path $P$, let $P = xuvy$. We claim that $u$ is a cut-vertex of $G$. Suppose this is not the case and let $Q$ be a shortest $x,y$-path in $G-u$. Since $\diam(G) = 3$ and $Q$ does not contain the edge $uv$, its length is at least $4$. Consider the subpath $Q'$ of $Q$ induced by its last four vertices, say $y_1$, $y_2$, $y_3$, and $y_4 = y$. It is possible that $y_3 = v$. We claim that $Q'$ is a shortest path also in $G$. Since we have assumed that $Q'$ is shortest in $G-u$, the only way $Q'$ is not shortest in $G$ would be that there exists a shortest path in $G$ between $y_1$ and $y_4=y$ using the vertex $u$. But this would mean that $uy\in E(G)$, which is not possible as $P$ is a diametral path. Hence $Q'$ is a shortest path in $G$ of length $3$ which does not contain $uv$. This contradiction proves the claim, that is, $u$ is a cut-vertex. By symmetry of $P$ we also deduce that $v$ is a cut-vertex.

Consider now the block $B$ of $G$ containing the vertices $u$ and $v$. Let $w\ne u,v$ be an arbitrary vertex from $B$. We first show that at least one of the edges $wu$ and
$wv$ exist. Suppose not. Then assume without loss of generality that $d_G(w,u) \le d_G(w,v)$. If $P$ is a shortest $w,u$-path, then this path (which must be of length $2$) together with an arbitrary edge between $u$ and a vertex not in $B$ is a diametral path which does not pass $uv$, a contradiction. Suppose next that $wu\in E(G)$ but $wv\notin E(G)$. Since $B$ is a block and $wv\notin E(G)$, there is another $w,v$-path in $B$ which does not pass $u$. Assuming that this is a shortest possible such path, the last two edges of it together with the edge between $v$ and a vertex adjacent to $v$ not from $B$ again gives a diametral path not containing $uv$. It follows that each of $u$ and $v$ is adjacent to every vertex of $B$. Then $B\setminus \{u,v\}$ is the graph $G_0$ required to show that $G\in {\cal G}_2$. Moreover, let $G_1$ be the graph induced by the vertices of $G$ which are not in $B$ and are closer to $u$ than to $v$. Then $u$ must be adjacent to each vertex of $G_1$ because otherwise we would have $\diam(G) \ge 4$. Analogously we can define $G_2$ with respect to $v$ and conclude that $v$ must be adjacent to each vertex of $G_2$. We conclude that $G\in {\cal G}_2$.
\qed

In the second part of this section we consider graphs with small edge general position number.

From Proposition~\ref{pro:earlier-results}(ii) we get that if $T$ is a tree with $m(T)\ge 2$, then $\gpe(T) = 2$ if and only if $T$ is a path. Moreover, if $G$ is not a tree, then considering a shortest cycle in $G$ which is necessarily isometric, we infer that $\gpe(G)\ge 3$. Hence we have the following easy result.

\begin{proposition}
\label{prop:gpe=2}
If $G$ is a graph with $m(G)\ge 2$, then $\gpe(G) = 2$ if and only if $G$ is a path.
\end{proposition}

We further have:

\begin{proposition}
\label{prop:gpe=3}
Let $G$ be a graph. Then $\gpe(G) = 3$ if and only if $G$ is $K_3$ or a tree with three leaves.
\end{proposition}

\proof
Clearly, $\gpe(K_3) = 3$. If $G$ is a tree with three leaves, then $\gpe(G) = 3$ by Proposition~\ref{pro:earlier-results}(ii).

Conversely, assume that $\gpe(G) = 3$. If $G$ is a tree, then the conclusion follows from Proposition~\ref{pro:earlier-results}(ii). Hence assume in the rest that $G$ is not a tree. Let $C$ be a shortest cycle of $G$. If $n(C)\ge 4$, then since $C$ is isometric, Proposition~\ref{pro:earlier-results}(i) implies that $\gpe(G)\ge 4$. Hence $C$ must be a triangle. Let $V(C)=\{x,y,z\}$. If $n(G)\geq 4$, then we may without loss of generality assume that $\deg_G(x)\geq 3$. Let $x'$ be adjacent to $x$, where $x'\not\in\{y,z\}$. Then $\{xy,yz,xz,xx'\}$ is an edge general position set of $G$. We conclude that $G = K_3$.
\qed

For graphs $G$ with $\gpe(G) = 4$ we have the following.

\begin{proposition}
\label{thm:necessary-gpe=4}
If $G$ is a graph with $\gpe(G) = 4$, then the following holds.
\begin{enumerate}
\item[(i)] $\Delta(G) \le 4$.
\item[(ii)] If $\Delta(G) = 4$, then $G$ is bipartite.
\end{enumerate}
\end{proposition}

\proof
We know from~\eqref{eq:Delta} that (i) holds. To prove (ii), consider an arbitrary vertex $v$ of $G$ with $\deg_G(v) = 4$ and consider the distance levels with respect to $v$. If there exists two adjacent vertices  in some of these distance levels, then the edge between them together with the edges incident with $v$ form an edge general position set. So we would have $\gpe(G)\ge 5$. We conclude that the distance levels of $v$ induce independent sets which in turn means that $G$ is bipartite.
\qed

Let $G_k$, $k\ge 1$, be the graph consisting of a chain of $k$ cycles $C_4$ sharing a vertex, see Fig.~\ref{fig:gpe=4} where $G_5$ is shown from which the formal definition of $G_k$ should be clear. Then it can be checked that $\gpe(G_k) = 4$ for each $k\ge 1$. More generally, the edge general position number remains $4$ if we replace each $C_4$ of $G_k$ by an arbitrary even cycle.

\begin{figure}[ht!]
\begin{center}
\begin{tikzpicture}[scale=1.5,style=thick]
\tikzstyle{every node}=[draw=none,fill=none]
\def\vr{3pt} 

\begin{scope}[yshift = 0cm, xshift = 0cm]
    \node [below=0.5mm] at (0,0) {};
    \node [below=0.5mm] at (1,0) {};
    \node [above=0.5mm] at (2,0) {};
    \node [below=0.5mm] at (3,0) {};
    \node [above=0.5mm] at (4,0){};
    \node [below=0.5mm] at (5,0) {};
    \node [below=0.5mm] at (0.5,0.5) {};
    \node [below=0.5mm] at (0.5,-0.5) {};
    \node [below=0.5mm] at (1.5,0.5) {};
    \node [below=0.5mm] at (1.5,-0.5) {};
    \node [below=0.5mm] at (3.5,0.5) {};
    \node [below=0.5mm] at (3.5,-0.5) {};
    \node [below=0.5mm] at (4.5,0.5) {};
    \node [below=0.5mm] at (4.5,-0.5) {};
    \node [below=0.5mm] at (2.5,0.5) {};
    \node [below=0.5mm] at (2.5,-0.5) {};

\path (0,0) coordinate (x1);
\path (0.5,0.5) coordinate (x2);
\path (1,0) coordinate (x3);
\path (1.5,0.5) coordinate (x4);
\path (2,0) coordinate (x5);
\path (3,0) coordinate (x6);
\path (3.5,0.5) coordinate (x7);
\path (4,0) coordinate (x8);
\path (4.5,0.5) coordinate (x9);
\path (5,0) coordinate (x10);
\path (4.5,-0.5) coordinate (x11);
\path (3.5,-0.5) coordinate (x12);
\path (1.5,-0.5) coordinate (x13);
\path (0.5,-0.5) coordinate (x14);
\path (2.5,0.5) coordinate (x15);
\path (2.5,-0.5) coordinate (x16);

\draw (x1)--(x2)--(x3) -- (x4)--(x5);
\draw (x6)--(x7)--(x8)--(x9)--(x10)--(x11)--(x8)--(x12)--(x6);
\draw (x5)--(x13)-- (x3) -- (x14)--(x1);
\draw (x5)--(x15)--(x6)--(x16)--(x5);

\draw (x1)  [fill=white] circle (\vr);
\draw (x2)  [fill=white] circle (\vr);
\draw (x3)  [fill=white] circle (\vr);
\draw (x4)  [fill=white] circle (\vr);
\draw (x5)  [fill=white] circle (\vr);
\draw (x6)  [fill=white] circle (\vr);
\draw (x7)  [fill=white] circle (\vr);
\draw (x8)  [fill=white] circle (\vr);
\draw (x9)  [fill=white] circle (\vr);
\draw (x10)  [fill=white] circle (\vr);
\draw (x11)  [fill=white] circle (\vr);
\draw (x12)  [fill=white] circle (\vr);
\draw (x13)  [fill=white] circle (\vr);
\draw (x14)  [fill=white] circle (\vr);
\draw (x15)  [fill=white] circle (\vr);
\draw (x16)  [fill=white] circle (\vr);

\end{scope}
\end{tikzpicture}

\end{center}
\caption{Graph $G_5$ with $\gpe=4$.}
\label{fig:gpe=4}
\end{figure}

There are many other examples of graphs $G$ with $\gpe(G) = 4$. For instance, the graph $G$ which is obtained from the disjoint union of $C_n$, $n\ge 6$, and $P_m$, $m\ge 2$, by identifying a vertex of $C_n$ by a leaf of $P_m$, also has $\gpe(G) = 4$.

\section{Edge general position sets in block graphs}
\label{sec:block}

While investigating in the previous section the graphs with extremal edge general position number, it turned out that cut-vertices are ubiquitous. In this section we thus consider the edge general position number in block graphs. To state our results, some preparation is needed.

We say that a block $B$ of a graph $G$ is {\em thick} if $n(B)\ge 3$ and that $B$ is {\em pendant} if it contains exactly one cut-vertex of $G$.
A graph is a \textit{block graph} if every block of it is complete. (See~\cite{furman-2022, henning-2022, li-2022} for some current research on block graphs.)

A vertex of a graph is \textit{simplicial} if its neighbourhood induces a complete subgraph. The set of simplicial vertices in a graph $G$ will be denoted by $S(G)$ and the cardinality of $S(G)$ by $s(G)$.
A block $B$ of a block graph is \textit{simplicial} if $B$ has at least one simplicial vertex.
An edge of a block graph is a \textit{simplicial edge} if it is incident with at least one simplicial vertex. The set of simplicial edges of $G$ will be denoted by $S'(G)$ and its cardinality by $s'(G)$.

Let $G$ be a block graph and let $B_1,\ldots,B_k$ be its simplicial blocks. Set $b_i = n(B_i)$ and $s_i = |S(G) \cap V(B_i)|$ for $i\in[k]$. Since an edge of $B_i$ is not simplicial if and only if its both endvertices are not simplicial, we infer that
\begin{equation}
\label{eq:s'}
s'(G)=\sum\limits_{i=1}^{k}\left[\binom{b_i}{2}-\binom{b_i-s_i}{2}\right].
\end{equation}

Let $G$ be a graph and let $p_1,\ldots,p_k$ be the consecutive vertices of a path $P$ in $G$, where $k\geq 2$. We say that $P$ is an \textit{internal path} of $G$ if $\deg_G(p_1)\geq 3$, $\deg_G(p_k)\geq 3$, and $\deg_G(p_i) = 2$ for $2\le i\le k-1$. We say that $P$ is a \textit{pendant path} of $G$ if $\deg_G(p_1)\geq 3$, $\deg_G(p_k)= 1$, and $\deg_G(p_i) = 2$ for $2\le i\le k-1$.

Let $G$ be a block graph. Then the {\em reduction} $R(G)$ of $G$ is constructed from $G$ as follows. For every internal path $P$ of $G$, identify the end vertices of $P$ (and remove all the inner vertices of $P$). Moreover, replace every pendant path $P$ by a pendant edge attached to $R(G)$ in the same vertex as $P$ is attached to $G$. Note that if $G$ is a path, then $R(G) = G$.

\begin{theorem}
\label{the:gpe-smoothing graph}
If $G$ is block graph, then $\gpe(G)=\gpe(R(G))$.
\end{theorem}

\proof
It is clear that an edge general position set of $R(G)$ yields an edge general position set of $G$ of the same cardinality. It follows that $\gpe(G) \ge \gpe(R(G))$. In the rest we thus need to prove the reverse inequality.

If $G$ is a path, then $R(G) = G$ and there is nothing to prove. Assume in the following that $G$ is a block graph which is not a path. Then $\deg_G(v)\geq 3$ for some $v\in V(G)$. It follows that $\gpe(G)\geq 3$.
Let $X$ be an arbitrary $\gpe$-set of $G$.
Let ${\cal P}$ be the set consisting of all internal and pendant paths of $G$. If $P\in {\cal P}$, then $P$ is isometric subgraph of $G$ and hence $|X\cap E(P)|\leq 2$. We hence distinguish the following two cases.

\medskip\noindent
{\bf Case 1.} $|X\cap E(P)|\leq 1$ for any $P\in {\cal P}$.\\
It is obvious that $X$ is also an edge general position set of $R(G)$ if $|X\cap E(P)|=0$ for any path $P\in {\cal P}$. Next, assume that $|X\cap E(P')|=1$ for some path $P'\in {\cal P}$ and let $e\in X\cap E(P')$. If $P'$ is a pendant path, then $e$ can be replaced by the pendant edge of $R(G)$ corresponding to $P'$, which keeps the property of being in edge general position. If $P'$ is an internal path, then let $u$ and $v$ be the end-vertices of $P'$. Let $u'$ and $v'$ be the respective neigbors of $u$ and $v$ which do not lie on $P'$. Then at most one of the edges $uu'$ and $vv'$ lie in $X$. Assume without loss of generality that $uu'\notin X$. Then $(X\setminus \{e\}) \cup \{uu'\}$ is also an edge general position set. Repeating this procedure we end up with an edge general position set $X'$ of $G$ in which no edge from the paths from ${\cal P}$ lies in $X'$. Since this set is also an edge general position set of $R(G)$ we conclude that $\gpe(G) \le \gpe(R(G))$.

\medskip\noindent
{\bf Case 2.} $|X\cap E(P)|=2$ for some path $P\in {\cal P}$.\\
Let $e,e'\in X\cap E(P)$. Let $e=xy$ and $e'=x'y'$.
Since $G$ is not a path, there exists a simplicial vertex $v$ such that $v\notin V(P)$. Let $v'$ be the neighbor of $v$ such that the edges $e$, $e'$, and $vv'$ lie on a common shortest path. Since $e,e'\in X$, we clearly have $vv'\notin X$. Now the set $X' = (X\setminus \{e\}) \cup \{vv'\}$ is an edge general position set of $G$. If $vv'$ is a pendant edge of a pendant path $P'$ of $G$, then $E(P') \cap X = \emptyset$, for otherwise an edge from $E(P') \cap X$ would lie on a shortest path together with $e$ and $e'$. It follows that when $X'$ is constructed, the number of pendant paths which contain two edges from an edge general position set is reduced by one. The same holds if $P$ is an internal path. Hence after repeating this procedure for every path from ${\cal P}$ with two edges from $X$ we end up with an edge general position set $X''$ of $G$ such that $|X''\cap E(P)|\leq 1$ holds for any $P\in {\cal P}$. Now we can proceed as in Case 1.
\qed

Next we prove general bounds for block graphs.

\begin{theorem}
\label{th:lower bound block graph}
If $G$ is a block graph, then $ s'(G)\leq \gpe(G)\leq \binom{s(G)}{2}+1$.
Moreover, the bounds are sharp.
\end{theorem}

\proof
Let $c_i$ be the number of cut-vertices from $B_i$ in $G$ and let $E_i$ be the set of simplicial edges from $B_i$ for $i\in [k]$.
Set $X=\bigcup\limits_{i=1}^{k} E_i$. Recall that $s'(G)=\sum\limits_{i=1}^{k}\big[\binom{b_i}{2}-\binom{b_i-s_i}{2}\big]$. We claim that $X$ is an edge general position set of $G$.

Consider an arbitrary shortest $u,u'$-path $P$ from $G$ and let $u=u_0,\ldots,u_\ell=u'$ be its consecutive vertices. Clearly, $P$ does not contain three edges of $X$ if $\ell\le 2$. Hence we may assume that $\ell \geq 3$. Since $G$ is a block graph, the shortest path $P$ is the unique shortest $u,u'$-path. If $u\in S(G)$, then it follows that $|E(P)\cap X|\leq 2$. Hence $X$ is an edge general position set of $G$.
Similarly, we also get the same result if $u'\in S(G)$.
If $u,u'\not\in S(G)$, then $|E(P)\cap X|=0$ and $X$ also is an edge general position set of $G$.
In consequence, we conclude that $\gpe(G)\geq s'(G)$.

To prove the upper bound we use induction on $t$, the number of blocks of $G$. If $t=1$, then $G$ is a complete graph, then hence $\gpe(G)= \binom{s(G)}{2}$ and the assertion holds.

Let now $t\ge 2$ and let $B$ be a pendant block of $G$. Let $G'$ be the block graph obtained from $G$ by removing the block $B$ but keeping the corresponding cut-vertex. Let $R$ be an edge general position set of $G$ of order $\gpe(G)$ and let $R'$ be the restriction of $R$ to $G'$. By the induction hypothesis, $\gpe(G') \le \binom{s(G')}{2}+1$.
Let $v$ be the vertex connecting $G'$ with $B$ in $G$ and let $s(B)$ be the set of simplicial vertices of $B$ in $G$. Then $s(B) = n(B) - 1$.
We distinguish the following two cases.

\medskip\noindent
{\bf Case 1.} $v$ is not simplicial in $G'$.\\
In this case $s(G)=s(G')+s(B)$ and hence $s(G')=s(G)+1-n(B)$. Since $n(B)<s(G)$, we can now estimate as follows:
\begin{align*}
|R|& \leq |R'|+|E(B)|\\
&\leq \binom{s(G')}{2}+1+|E(B)|\\
&= \binom{s(G)+1-n(B)}{2}+1+\binom{n(B)}{2}\\
&=\frac{(s(G)+1-n(B))(s(G)-n(B))}{2}+\frac{n(B)(n(B)-1)}{2}+1\\
&= \binom{s(G)}{2}+1+(n(B)-s(G))(n(B)-1)\\
&\leq \binom{s(G)}{2}+1\,.
\end{align*}

\medskip\noindent
{\bf Case 2.} $v$ is simplicial in $G'$.\\
In this case $s(G)=s(G')-1+s(B)$, which means that $s(G')=s(G)-n(B)+2$.
If $n(B)\geq 3$, then
\begin{align*}
|R|& \leq |R'|+|E(B)|\\
&= \binom{s(G')}{2}+1+|E(B)|\\
&= \binom{s(G)+2-n(B)}{2}+1+\binom{n(B)}{2}\\
&= \binom{s(G)}{2} + 1 + (s(G)-n(B)) (2-n(B)) + 1\,.
\end{align*}
Since $n(B)\geq 3$, we have $2-n(B)<0$. Moreover, because $s(G)-n(B)\ge 0$, we can continue the above estimation with $|R| < \binom{s(G)}{2} + 2 - 1$, hence $|R| \le  \binom{s(G)}{2} + 1$.

If $n(B)= 2$, it follows that $s(G)=s(G')$.
Assume, without loss of generality, that $V(B)=\{v,v'\}$.
If further $vv'\not\in R$, then it is obvious that $|R|=|R'|\leq \binom{s(G)}{2}+1$.
Hence we may assume that $vv'\in R$.
If $k=2$, we see that $s(G)=n(G)-1$ and we can conclude that $|R|\leq \binom{s(G)}{2}+1$.
If $k\geq 3$, let $v''$ be a neighbor of $v$ in $G'$.
Then there exists one shortest $v,x$-path $P$ in $G'$ such that $|E(P)\cap R'|=2$.
Otherwise $R'\cup \{vv''\}$ is an edge general position set of $G'$ contradicting the maximality of $R'$ in $G'$.
Hence we have $|R|\leq |R'|\leq \binom{s(G)}{2}+1$, completing the proof of the upper bound.

The lower bound is sharp on trees. Indeed, if $T$ is a tree, then $\gpe(T)$ is the number of leaves of $T$, and this is equal to $s'(T)$. To show that the upper bound is sharp, consider a graph $G$ obtained from the complete graph $K_n$, $n\ge 2$, by attaching selecting $k\in [n]$ vertices of $K_n$ and add a separate pendant path (of arbitrary length) to each of the selected vertices. (See Fig.~\ref{fig:G'-and-G''}, where the graph $G''$ is such a graph.) Then we infer that $s(G) = n$. Moreover, the edges of $K_n$ together with an arbitrary additional edge of $G$ form a largest edge general position set of $G$. Indeed, this set is clearly an edge general position set. To see that it is largest, we can argue that if an edge general position set $X$ of $G$ contains at least two edges from $E(G)\setminus E(K_n)$, then all but one of these edges can be replaced by a unique edge from $E(K_n)\setminus X$ keeping the property of being an edge general position set.
\qed

\begin{figure}[ht!]
\begin{center}
\begin{tikzpicture}[scale=1.3,style=thick]
\tikzstyle{every node}=[draw=none,fill=none]
\def\vr{3pt} 

\begin{scope}[yshift = 0cm, xshift = 0cm]
    \node [below=0.5mm] at (1,0) {};
    \node [below=0.5mm] at (2,0) {};
    \node [above=0.5mm] at (3,0) {};
    \node [below=0.5mm] at (3.75,0.5) {};
    \node [above=0.5mm] at (3.75,-0.5){};
    \node [below=0.5mm] at (1.5,0.5) {};
    \node [below=0.5mm] at (1.5,-0.5) {};

\path (1,0) coordinate (x2);
\path (2,0) coordinate (x3);
\path (3,0) coordinate (x4);
\path (3.75,0.5) coordinate (x5);
\path (3.75,-0.5) coordinate (x6);
\path (1.5,0.5) coordinate (x7);
\path (1.5,-0.5) coordinate (x8);

\draw (x2)--(x3)--(x4) -- (x5)--(x6)--(x4);
\draw (x7)--(x8);
\draw (x2)--(x7)-- (x3) -- (x8)--(x2);

\draw (x2)  [fill=white] circle (\vr);
\draw (x3)  [fill=white] circle (\vr);
\draw (x4)  [fill=white] circle (\vr);
\draw (x5)  [fill=white] circle (\vr);
\draw (x6)  [fill=white] circle (\vr);
\draw (x7)  [fill=white] circle (\vr);
\draw (x8)  [fill=white] circle (\vr);

\draw (2.5,-1) node {$G'$};
\end{scope}
\end{tikzpicture}
\hspace{10mm}
\begin{tikzpicture}[scale=1.35,style=thick]
\tikzstyle{every node}=[draw=none,fill=none]
\def\vr{3pt} 

\begin{scope}[yshift = 0cm, xshift = 0cm]
 \node [below=0.5mm] at (0,0) {};
    \node [below=0.5mm] at (1,0) {};
    \node [below=0.5mm] at (2,0) {};
    \node [above=0.5mm] at (3,0) {};
    \node [below=0.5mm] at (4,0) {};
    \node [below=0.5mm] at (1.5,0.5) {};
    \node [below=0.5mm] at (1.5,-0.5) {};

\path (0,0) coordinate (x1);
\path (1,0) coordinate (x2);
\path (2,0) coordinate (x3);
\path (3,0) coordinate (x4);
\path (4,0) coordinate (x5);
\path (1.5,0.5) coordinate (x6);
\path (1.5,-0.5) coordinate (x7);
\draw (x1)-- (x2)--(x3)--(x4) -- (x5);
\draw (x2)--(x6)-- (x7);
\draw (x6)--(x3)-- (x7);
\draw (x2)--(x7);

\draw (x1)  [fill=white] circle (\vr);
\draw (x2)  [fill=white] circle (\vr);
\draw (x3)  [fill=white] circle (\vr);
\draw (x4)  [fill=white] circle (\vr);
\draw (x5)  [fill=white] circle (\vr);
\draw (x6)  [fill=white] circle (\vr);
\draw (x7)  [fill=white] circle (\vr);

\draw (2.5,-1) node {$G''$};
\end{scope}
\end{tikzpicture}
\end{center}
\caption{Block graphs $G'$ and $G''$ with $\gpe(G') = 9$ and $\gpe(G'') = 7$.}
\label{fig:G'-and-G''}
\end{figure}

We follow with a large class of block graphs which attain the lower bound of Theorem~\ref{th:lower bound block graph}. It is defined as follows. We say that a graph $G$ is a \textit{thick-leaved tree} if $G$ is obtained from a tree $T$ by a sequence of the following operations. Let $v$ be a leaf of $T$ and let $v'$ be its unique neighbor. Then replace the vertex $v$ by a complete graph $K$, and replace the edge $v'v$ by an edge between $v'$ and one vertex of $K$. See Fig.~\ref{fig:G'-and-G''} where $G'$ is a thick-leaved tree. Note that a thick-leaved tree is a block graph. Note further that a pendant block of a thick-leaved tree is either a thick block or a $K_2$. Moreover, its simplicial blocks coincide with its pendant blocks.

\begin{proposition}
\label{prop:lower bound}
If $G$ is a thick-leaved tree and $B_1,\ldots,B_k$ are its simplicial blocks, then
$$\gpe(G)=\sum\limits_{i=1}^{k}\binom{n(B_i)}{2}.$$
\end{proposition}

\proof
By Theorem~\ref{th:lower bound block graph}, we know that $\gpe(G) \ge s'(G)$. Since $s'(G) = \sum\limits_{i=1}^{k}\binom{n(B_i)}{2}$, the lower bound follows.

To prove that $\gpe(G) \le \sum\limits_{i=1}^{k}\binom{n(B_i)}{2}$, consider an arbitrary edge general position set $X$ of $G$ and suppose $e=uv\in X$, where $u$ and $v$ are cut-vertices of $G$. Let $G_u$ and $G_v$ be the two components of $G-e$, where $u\in G_u$ and $v\in G_v$. Then it follows that $X\cap E(G_u)=\emptyset$ or $X\cap E(G_v)=\emptyset$. Indeed, otherwise there exist edges $e'\in X\cap E(G_u)$ and $e''\in X\cap E(G_v)$, but then $e'$, $e$, and $e''$ would lie on the same shortest path in $G$.
Hence we may assume, without loss of generality, that $X\cap E(G_u)=\emptyset$.
Let $f$ be an arbitrary simplicial edge of $G_u$. Note that $f$ may be adjacent to $e$.
Then it is clear that $(X\setminus\{e\})\cup \{f\}$ also is an edge general position set of $G$. Repeating this process we end up with an edge general position set of $G$ which contains only simplicial edges and has the same cardinality as $X$. We conclude that $|X| \le s'(G)$.
\qed

\section*{Acknowledgements}

This work has been supported by T\"{U}B\.{I}TAK and the Slovenian Research Agency under grant numbers 122N184 and BI-TR/22-24-20, respectively. Sandi Klav\v{z}ar also acknowledges the financial support from the Slovenian Research Agency (research core funding P1-0297 and projects J1-2452 and N1-0285). Jing Tian has been supported by NNSF of China (Grant No.\ 12271251).

\section*{Declaration of interests}

The authors declare that they have no known competing financial interests or personal relationships that could have appeared to influence the work reported in this paper.

\section*{Data availability}

Our manuscript has no associated data.


\end{document}